# On the difficulties of acquiring mathematical experience: Case rural education

*Em Português:* **Sobre as dificuldades na aquisição da experiência matemática: Caso da educação rural**

by Bernhelm Booß-Bavnbek (Roskilde University, Denmark)[1]

*Abstract.* Based on a variety of philosophical approaches and my own work for decades in pure and applied mathematics teaching and research, I explain my view upon the basic difficulties of acquiring the "Mathematical Experience" (in the sense of P.J. Davis and R. Hersh 1981) and submit a list of claims how these difficulties can and should be confronted.

*Em Português Resumo*: Com base em diferentes abordagens teóricas e em meu próprio trabalho de ensino e pesquisa, que venho desenvolvendo por décadas, em matemática pura e aplicada, explico meu ponto de vista sobre as dificuldades básicas de aquisição da "Experiência Matemática" (no sentido atribuído por P.J. Davis and R. Hersh 1981) e apresento uma lista de afirmações sobre como essas dificuldades podem e devem ser confrontadas.

**Introduction – Bitterness and seriousness of the topic "Rural education and mathematics".** On the one side, the rural environment is the oldest *form of life* (Antonio Gramsci, 1891-1937, in his Prison Notebooks *Quaderni del carcere* of 1929-1935 in reference to Wittgenstein) of human beings. Living in the countryside or in marginalized suburban regions is still the dominant form of life for a substantial and often underprivileged part of the world population, screaming for a *humanist and libertarian pedagogy* (Freire 1972, p. 40, see also the moving recent elaboration Mesquita 2013 in this issue). On the other side, numbers, calculations, geometry, mathematics have been the precursors of all sciences and technology (Damerow and Lefèvre 1981); being visible - or hidden in black boxes and chipified (Keitel, Kotzmann & Skovsmose 1993), mathematics plays an ever greater role in designing products and production processes, formatting social relations, and puffing-up all modern research. It is by far the oldest teaching subject (Høyrup 2002) and the largest through history.

In view of that, rural math teachers are confronted with the Kantian *categorical imperative* in Karl Marx' reformulation "to overthrow all those conditions in which man is an abased, enslaved, abandoned, contemptible being" (1844). How can and must mathematics education contribute to that humanistic agenda of rural education and what are the greatest obstructions to overcome on that path?

---

[1] Email: booss@ruc.dk, mail: NSM/IMFUFA, Roskilde Universitet, Postboks 260, DK-4000 Roskilde, Denmark

## Approach 1 – The simplicity of mathematical experience and its natural character.

People from rural areas, pupils and their parents deserve respect. They should no longer be addressed as the feeble, fail and shaky. Teachers must learn to appreciate and appraise the value of practical knowledge and to deal with the highly complex environment of rural society that fosters two contrary positions: Possibly, most students in a peripheral Brazilian school will agree with Pink Floyd's "We don't need no education. We don't need no thought control. No dark sarcasm in the classroom. Teachers leave them kids alone". These pupils will be blocked by poverty, extensive child labor and the foreignness of teacher and classroom set-up. Most probably, they will not believe in schools. They will not care about contents because they know what is waiting for them at the end of schools: unemployment or labor exploitation. Other pupils with more resources will take another, more modern rejection position related to the ubiquity of the internet and proudness based on the competent daily handling of the technologically advanced machinery operated in the countryside. That can give pupils and parents the feeling that further school training is dispensable.

### Approach 1 example
The world famous mathematician Izrail M. Gelfand (1913–2009) insisted upon the basic simplicity of mathematics. Gelfand never thought that there was anything very difficult about math, regarding it as "a way of thinking in everyday life". "You can explain fractions even to heavy drinkers," he told an interviewer in 2003. "If you ask them, 'which is larger, 2/3 or 3/5?' it is likely they will not know. But if you ask, 'Which is better, two bottles of vodka for three people, or three bottles of vodka for five people?' they will answer you immediately. They will say two for three, of course." (From the obituary in the Daily Telegraph of 26 Oct 2009).
   Heavy drinkers aside, from ethnomathematical studies there is ample evidence about the deep roots of arithmetic and geometry in people's trade, craft and art and their remarkable achievements (D'Ambrosio 2008, Gerdes 2000, Knijnik and Wanderer2010).

### Approach 1 theoretical explanation
One may rightly doubt Gelfand's claim that the way of thinking in everyday life is the basis of scientific thinking. For me and many other mathematicians, his way of thinking remains a mystery and a treasure with its unlimited fecundity and its surprising revelations of exotic interconnections. Ordinary mathematicians, not to speak of ordinary people, will hardly recognize any traces of common sense and everyday life thinking in Gelfand's work.  Below, I shall come back to the sharp distinction between everyday language and mathematical formalism.
As long as we do not exaggerate Gelfand's claim, there is a theoretical explanation though, for the basic mathematical giftedness of all human beings, at least if we are willing to follow Noam Chomsky's lifelong investigation of the innate generative grammar underlying human language acquisition (Pinker 1994). Consequently, we all have passed our hardest math test in our life at the age of 2 or 3, when we became able to analyze and reproduce the complicated algebraic patterns of the native language.

### Approach 1 practical conclusion and claim
Recall the apocryphal statement attributed to Napoléon Bonaparte (1769-1821): "Tout soldat français porte dans sa giberne le bâton de maréchal de France" (Every French soldier carries in his knapsack the baton of a Marshal of France)! As the historians tell us, the key element in the Napoleonic reforms was a greater integration of the populace into the mechanisms of the state: at that time, its army, and its administration in general. This development let air into the stymied

French social system. Mathematics and, equally well, technological education can and has to do a corresponding job under present conditions, namely to encourage and integrate pupils in rural areas.

## Approach 2 – The difficulty of mathematical experience and its formal and abstract character.

Mathematics is formal and abstract. That constitutes the power of mathematics. Shortly put, abstract formulas serve efficiently as the memory of mathematically encoded human experiences. At the same time, mathematical formulas provide a field of imagination and permit transfer, variation and adaption of experiences from one field of human activity to another. That's the plus-side. On the negative side, math teachers, pupils and the general public all agree that mathematics is difficult to learn; new mathematical concepts are difficult to grasp; formulas and graphs are difficult to read; and abstract arguments are difficult to follow.

### Approach 2 example

Normal human beings need time and concentration as long as we don't understand the meaning of the abstract signs and arguments. We easily get confused and sometimes stressed. Later, however, when the penny has dropped, everything becomes incredibly easy and the power of abstractions and formulas may become evident. Then teacher and pupil will share the mathematicians' ultimate joy: "There is no feeling quite like that which comes after you have proved a good theorem, or solved a problem that you have worked on for a long time. Driven by the heat of passion, the words burst forth from your pen, the definitions get punched into shape, the proofs are built and bent and patched and shored up…" (Krantz 1997, p. xi). Then we meet a new challenge, namely how to understand our now surpassed previous difficulty. The young rural math learner's multifaceted submerging process mirrors the mixed feelings of even the most experienced and ingenious professional mathematician when reading and comprehending a new mathematical paper as referee for evaluation: *Incomprehensible! – No, wrong! – No, trivial (I did it before)!* I could quote mathematical heroes like Carl Friedrich Gauss (1777-1855) or Lars Hörmander (1931-2012).

### Approach 2 theoretical explanation

Since Aristotle (384 – 322 BCE), logicians have been convinced that concepts and arguments can be clarified by wrapping them into formal abstract expressions. Many of them may agree with J. Barkley Rosser's (1953, p. 7) famous – and to me a bit exaggerated - claim, that "once the proof is discovered, and stated in symbolic logic, it can be checked by a moron." In social sciences it is even common praxis to attribute enhanced credibility to one's own arguments by dressing them up with diagrams and abstract formulas. It seems to me that on the contrary most mathematicians would rather support Yuri Manin's (2010, p. 36) laconic statement: "The human mind is not at all well suited for analyzing formal *texts*." Note that Skovsmose (1994, p. 48) states no less laconically, when discussing the *formatting power* of mathematics upon society and what he calls the *Vico paradox*, that "we seem completely unable to establish such an understanding (sc. of our artefacts) in the case of technology. We seem to be without the capacity of grasping the limits and the full consequences of our technological enterprises… Although technology is a human construction, we do not seem to possess the capacity to comprehend what we have constructed."

Manin (2010, p. 33) explains his hypothesis on several levels: He compares logic and mathematical formalism with rigid stencils which we can impose on any artificial system like mathematics, electrical network designs or computer programs, but not on reality. As we know from physics, reality requires the application of a variety of partly overlapping and mutually not always consistent approaches. "The physicist's descriptions do not have to form a consistent or coherent whole; his job is to describe nature effectively on certain levels. Natural languages and the spontaneous workings of the mind are even less logical. In general, adherence to logical principles is only a condition for effectiveness in certain narrowly specialized spheres of human endeavor." Manin elaborates seemingly parallels between expressions of mathematical logic and everyday language and shows how misleading such parallels are. He concedes, though, that formal mathematics, in which a single contradiction destroys the entire system, clearly has the features of poetic hyperbole. Similarly, but more generally, Hans Freudenthal (2002, p.2) saw "mathematics as an art, a mental art to be sure, which for most people will be closer to crafts than to sciences, a tool rather than an aim in itself, more relevant because it works than because it is certain. ... Although many people trust mathematics more than it deserves, it works only when it is rightly applied.... Once one has admitted that mathematics is an art, one cannot shirk the responsibility of judging whether, in particular cases, it is being properly used or rather being abused."

Repeatedly, Manin (2010, p. 34) emphasizes that "the choice of the primitive modes of expression in the logic of predicates does not reflect psychological reality. *Elementary* logical operations, even one-step deductions, may require a highly trained intellect; yet, logically *complicated* operations can often be performed as a single elementary act of thought even by a damaged brain." Perhaps most strikingly, Manin (l.c., pp. 34-36) illustrates his view by presenting details of the study of a patient with brain injuries, carried out by the psychologist A. R. Luria (1972).

I have no reason to doubt Chomsky's findings that the capacity to make abstractions is deeply rooted in our genes, as explained above. That does not contradict Luria's and Manin's finding that abstractions typically are difficult to grasp and may even be conflicting with everyday perception.

## Approach 2 practical conclusion and claim

In a particularly extreme way, people from rural areas are exposed to the cultural clash immanent in abstractions, formalism and symbol processing. Teachers must help them to experience that clash as a positive step like processes of adolescence and not as a series of defeats. However, it doesn't help with well-intended lies or self-deception about easy access to mathematical abstractions. Acquiring mathematical experience is nothing that falls from heaven or comes from playing on the ground. It requires work, concentration, exercises, and endurance: Ὁ μὴ δαρεὶς ἄνθρωπος οὐ παιδεύεται (The non-flayed human will not be educated, Menander, c. 341/42– c. 290 BCE, disseminated by J.W. Goethe as motto over his autobiography *Dichtung und Wahrheit*), or less draconic, *Ohne Fleiss kein Preis* (Without hardworking no praise, after Hesiod, thought by scholars to have been active between 750 and 650 BCE).

## Approach 3 – The semiotic and existentialist perspective: Mathematical experience as language and as notation and the teacher's role.

In Approach 1, I emphasized the intimate relation between mathematical concepts and everyday knowledge. A keyword was *language*: nursing perceptions of another human being by language games, like moving around in different old towns (Wittgenstein 1953 in his lifelong endeavor to correct the logistic aberrations of his young days). In Approach 2, I emphasized the dual phenomenon, i.e., the distinct character of mathematical abstractions and the sharp difference between formal arguments, computer languages and programs on the one side and common language on the opposite side. Here the keywords were *abstraction* and *notation*. In the present section I aim to mediate between the two complementary concepts.

### Approach 3 example

Again and again, clever people have tried to wipe out the distinction. Perhaps the most baroque try is related to the concept of *Artificial Intelligence* (*AI*). Rightly, Minsky (1967) may be considered as AI's founding manifesto, in particular in its unconfined willingness to promise everything funding agencies may wish, quote: "within a generation ... the problem of creating 'artificial intelligence' will substantially be solved" (l.c., p. 2). That propaganda worked indeed for gathering substantial funding for years to come, mostly from the US Department of Defense. However, nothing of the early promises was delivered in the 50 years since then: no automatic language translation, as [Google translate] can witness, and no automatic pattern recognition as witnessed by the inverse Turing test *Completely Automated Public Turing test to tell Computers and Humans Apart* (*CAPTCHA*), widely applied to prevent hacking of internet pages. With hindsight we can see that AI's total failure in its core promises can be explained by underestimating the differences between our approaches 1 and 2. As long as AI stays solely with symbol processing, as indicated in our Approach 2, impressive achievements are attainable: data mining and the management of huge knowledge data bases; flexible programmable robots; programs that beat chess champions. Computers can even be programmed to pass the Turing test in Turing's original and sarcastic meaning, namely simulating a British Lordship's understanding and speaking so well that the Lordship himself would never guess that there is a machine and not another lordship at the other end of the line. For its naïve mix of our Aspects 1 and 2, for its lack of psychology and lack of neurological insight (not obtainable then and perhaps not available now), AI's core program had to go down.

### Approach 3 theoretical explanation

What then could be a more viable reconciliation between our two Aspects, and hopefully a productive one for teaching mathematical experiences in the countryside? In this Note, I shall argue for the semiotic view due to the US-American logician, mathematician and philosopher Charles Sanders Peirce (1839 –1914). In his huge work, Peirce summarizes classical European philosophy with special interest for scholastic philosophers, Scottish realism and Kant (Peirce 1877 in my free edition may give a taste of his way of arguing).

Roughly speaking, Peirce considers both reality and our thinking, speaking and writing as sign systems. The two sign systems are extremely different. It needs an observer to relate *the sign systems of reality* (our senses' impressions and measurements) with *the formal sign systems* of our pre-knowledge. The quality of the relation cannot be determined by maximal coincidence. On the contrary, only the distance of the formal signs from the observed signs makes formal signs valuable for practical goals. So, the usefulness for *practical goals of a specific observer* or thinker is the only valuable criterion for the adequacy and fecundity of any relation between formal and real signs. That puts the observer, i.e., the mediator between the sign systems, after all the teacher and the learner in the center of any epistemological and learning process.

Here the next philosopher comes into play, the Danish existentialist and theologian Søren Kierkegaard (1813 – 1855). We may adapt his analysis of Mozart's Don Giovanni in *Either/or* to the classroom situation. Equally well as, according to Kierkegaard, the Christian God is depicted in the seducer Don Juan, we may also identify the teacher's role with Don Juan's seductive arts towards Zerlina, as soon as he was alone with her in the Duet: "Là ci darem la mano" – "There we will entwine our hands". Contrary to that, the two other main male figures of the opera, Don Juan's servant Leporello and Donna Anna's lover Don Ottavio don't show that seducing quality and interpersonal intensity necessary for successful teaching. Leporello is good in recalling and summarizing all of Don Juan's adventures: "Madamina, il catalogo è questo" – "My dear lady, this is the catalogue". He offers numbers instead of passion. Don Ottavio is great in swearing love and revenge: "Ah, vendicar, se il puoi, giura quel sangue ognor!" – "Ah, swear to avenge that blood if you can!" But then nothing happens. Only Don Juan's existential full presence and his acceptance of the risk of rejection can provide a role model for the teacher in combining language games and formalism for and with the pupils.

An important and enormously inspiring modern follower of Peirce (via the pragmatist philosophers William James and Harald Høffding and the physicist Niels Bohr) is the Danish multimath, astronomer, computer scientist and psychologist Peter Naur. Like Wittgenstein, Naur combines our listed two aspects in his own person, biography and work: On one side, his last name is the N in the BNF notation (Backus-Naur form), used in the description of the syntax for most programming languages and of some similarity with Chomsky's generative grammar. He was among the creators of ALGOL which became the model for many modern computer programming notations. On the other side, the same Naur (1992) used his last 30-40 years to address the psychological aspects of software design; the interplay with the users; the programmers' prejudices and the underlying perception of the treated reality.

Within the math education community, it seems to me that Otte (1994, 1997) comes close to the Peircean view, while Radford (2006) provides an interesting (yet to me not fully convincing) critical re-examination of Peirce's and Husserl's (1970) attempts at reconciling the subjective dimension of knowing with the alleged transcendental nature of mathematical objects.

## Approach 3 practical conclusion and claim

The teacher's intensity must become felt and visible in the classroom for a credible combination of language game and formalism. Pupils and parents get many bits of information from radio and TV which will be forgotten immediately. The teacher's role is to say something the pupils will remember for the rest of their lives. Mathematics must, via the teacher's personality, appeal to the feelings of the pupils. The surrounding society must support the teacher by conferring the corresponding required authority and trust. I quote the Finnish example (Sahlberg 2011, p. 130): "The culture of trust meant that education authorities and political leaders believe that teachers, together with principals, parents, and their communities, know how to provide the best possible education for their children and youth. Trust can only flourish in an environment that is built upon honesty, confidence, professionalism, and good governance" (after years of terrible civil war and stagnation, we may add). In the Finnish example only students with top grades can enter the teacher education; the teachers' professional autonomy in schools and their classrooms is respected; and their career choice not called into question. Outside inspectors are barred to judge the quality of their work, and no merit-based compensation policy influenced by external measures is imposed. (ibidem, p. 93-95).

Moreover, it seems to me that the semiotic challenges and the existentialist role of the math teacher define serious limitations to the replacement of expensive teachers by cheap machines for online interaction with the pupils. Though, admittedly, there are numerous well-documented positive experiences with preserving and even enhancing (mathematics teachers') reflections in online (teacher) training (Sánchez 2010).

## Approach 4 – The mathematical experience and the pragmatics of curriculum discussions.

Talking to math teachers of any age and on any level, most people will be at once impressed by their passion for their job and the richness of the individual teacher's experience: their special secrets how to move their pupils forward. To facilitate communicating teaching concepts and experiences among teachers and to pupils and students at the new founded Roskilde University, in 1972 three maxims for the learning of mathematics and sciences were suggested: *in, with, about*. See also (Skovsmose 1994) for a professional philosophical definition of the three maxims.

### Approach 4 example

When learning about plane figures, triangles, squares, angles and lengths, the *in-maxim* might emphasize the simplicity of basic concepts, the difference between symmetry and general position and the invariant validity of the Pythagorean Theorem which is trivial for isosceles rectangular triangles, as the philosopher Arthur Schopenhauer noticed (1788-1860) and erroneously perceived that as the content of the theorem. Here the topic is the joy of pure mathematical thinking, i.e., the joy of exploring the capacity of one's own brain in almost pure laboratory-type environment.

The *with-maxim* might emphasize the determination of unknown distances and heights in the landscape and other applications to experience the power of calculating with letters and geometric figures instead of numbers.

The *about-maxim* might invite to reflect about past history when triangle geometry was a key tool for military and civilian terrain measurements and the scientific, technological and societal changes that accompany their replacement by ready-to-use GPS and laser tools.

### Approach 4 theoretical explanation

For centuries, even millennia, math teachers could teach in peace without being exposed to changing didactics commands and fashions every few years. For pupils of the upper classes, the *in-maxim* was dominant: The future administrator had to learn and internalize the feeling of superiority that can be acquired from the knowledge of methods of solution for complicated mathematical exercises. That can be read from the Old Babylonian clay tablets of (Høyrup 2002) and is confirmed by math teaching in Greek antiquity. Math belonged to the pentagonal unity of possibility games, jointly with philosophy, travelling, theater and democracy. The *with-maxim* of the performance of practical calculations, e.g., for land property evaluation, military supplies or commerce was separated and left to the practitioner's schools.

Roughly speaking, the predominance of the *in-maxim* lasted long up to the end of the last century, though supplemented by punctual *with-maxims* regarding single topics like the introduction of stereometry and the discussion of differentiable curves with their minima and maxima and systematic *about-maxim* insisting on half-religious platonic adoration of the beauty and mysteries of math.

It should be noted that socialist writers repeatedly expressed their distrust in the predominance of the *in-maxim* (Lenin 1909). They perceived the in-maxim as socially elitist. Reactionary governments, on the other hand, were suspicious when math teachers opened the class room to real world problems (Schubring 2007).

Special attentions deserve the intricate paths to form new mathematical concepts in the brain of the pupils, foremost the *reification* (Sfard 1991). Consider, e.g., the formation of the concept of a "function": In their rural environment, pupils have seen various tables, graphs and formulas regarding crops, energy consumption, prices, wages, time, movements of a thrown stone or of the moon. Working with these different representations; adding or subtracting two graphs with the same independent variable and a comparable depending variable; scaling curves up and down; linearizing one and smoothing another to an exponential law, all these activities will make it natural for the student to give the class of objects a name – function; then to lift their work to a higher level; and then to slough the bindings to the various contexts. Working further will support inducing additional and for the pupils new concepts along the same path of reification. Such considerations, concepts and content may appear as an academic matter. They are not. They help students to discover new aspects of their own brain and their form of life.

For the last 100 years, math teachers, more than teachers of any other subject, have been exposed to radically changing educational maxims: In addition to the wrangling between the in- and with-maxim, there were the repeated attempts to balance geometry vs. the function concept; discrete vs. continuous math.; the conduct of a fixed syllabus vs. a pupil-oriented learning; the training of skills vs. the acquisition of competences.

I shall only comment on one controversy, the balance between structural foundations and mathematical substance. After the Sputnik shock and the impressive reconstruction of the war-devastated Russia and its allies without help from abroad there was a need felt in the West to catch up in mathematics and technology education with the Soviet Union in the shortest possible time and in the most effective way. In hindsight, it was natural to focus on generalities and structures when one felt to be under time pressure and confronted with immense challenges. While that orientation was controversial within the math community and shook the wider public of helpless parents, it may have been the best to do in those early years of computerizing and symbol processing. We will never know: The reform was halted before it could show possible merits (Moon 1986).

### Approach 4 practical conclusion and claim
In everyday teaching and learning, the categories of in, with and about may not always be clearly distinguishable. A clear declaration of the game played may help the pupils to find out whether capacities of their brain are discussed - or properties and design of technical and social inventions and constructions - or preconditions and consequences of mathematical procedures for humans and in society.

### Approach 5 – The necessity of reflection and contextual discourse for acquiring some mathematical experience.
Here I shall address a deep contradiction in all math teaching: the basically context-free character of mathematical formalisms and the necessity to embed them into concrete contexts for the purpose of learning and understanding their meaning. For a moment I shall leave the old controversy aside where mathematical concepts come from; to what extent they are abstractions from real world experiences and to what extent they are free inventions of the human mind, possibly guided by dreams of a platonic world. The basic contradiction is much more serious. As emphasized above, the power of mathematics comes from its versatility, from the potential of transfer of experiences from one human field to another one by abstractions, symbols, formalism. Thus, freeing mathematical formalisms from contexts is mandatory for making them applicable in new spheres of human activity. In this aspect, my position is orthogonal to some recent approaches in math education (sketched, e.g., in Miguel and Mendes 2010): Inspired by the French philosophers J. Derrida and M. Foucault, the authors emphasize the essential difference between the various language games and contexts. Addressing contexts and differences between contexts should make the main parts of the math teaching. Of course they are right in many concrete teaching situations.

However, once again I shall emphasize that math is distinguished by imposing the same or similar rules to often very different looking contexts. For the mathematician, the equations for planetary motion and a ball rolling down a skew plane are the same; the ups and downs of the spread of a contagious disease are analyzed by models derived from the pendulum oscillations; chaotic behavior is generated by simple 3d degree complex polynomial; the basic mode for medical tests are casts of a dice. Where *postmodernists* emphasize differences of contexts, the mathematician sees basic similarities. In so far, math is similar to religion, the other reality oriented and unity creating or claiming ideology. Compared to religion, math has, however, the advantage that its application in real situations after all coincides with our experiences, while I, personally, cannot confirm the coincidence of religious experience with reality. As Philip Davis and I explained at another place recently, one should, however, not exaggerate the claim of unity of and via mathematics (Booß-Bavnbek and Davis 2013).

Clearly, and in full concordance with (Miguel and Mendes 2010) and (Radford 2006), the meaning of abstractions will become accessible, though, only through contexts, and finally through action in specific contexts. That was what Karl Marx (1818-1883), the *greatest philosopher of pragmatism* according to B. Russell, wrote in the famous $11^{th}$ Feuerbach thesis (Marx 1845): "Die Philosophen haben die Welt nur verschieden interpretiert; es kommt darauf an, sie zu verändern." (Philosophers have hitherto only interpreted the world in various ways; the point is to change it.) By the way, it is easy to read from Marx' iconic, often reproduced handwriting of the thesis, that in the German original there is no "aber" (no "though") after "kommt", contrary to the faked and false quote on http://en.wikipedia.org/wiki/Theses_on_Feuerbach, Footnote 5. In teaching mathematical experiences (as well as in political praxis), we should not place the second part of Marx' Feuerbach thesis over its first part. We should discard all well-intended but uninformed actionism and contextualization. We have to insist that there is no meaning without content and method (like in the beautiful title and content of Alexandrov et al. 1963) and no meaningful action without thorough investigation and understanding in a broader setting beyond an acute context.

## Approach 5 example 1
Math teachers handle this contradiction in different ways. They may neglect relevant contexts like weather prediction, industrial design, telecommunication, computer architecture, construction and civil engineering and other fields that have their own subject tradition and often are too complex for treatment in class. Instead of that, typically less representative contexts are chosen that can have the advantage of easier accessibility. A teacher may, e.g., judge that treating the gravitational laws of planetary movement is mathematically too demanding for a class, and prefer to discuss mathematically similar relations in respect to the attractiveness of urban shopping centers. Indeed, some of the math applied can be similar. However, the pupils' impression of the character of mathematical modeling will be blurred by such pedagogically well-intended compromises: the decisive difference between ad-hoc assumptions and models based on first principles will be wiped out.

## Approach 5 example 2
Teachers also go astray when they identify themselves too much with a mathematically nicely formulated modeling task from faked real world which lacks realism. Attentive pupils will feel fooled. Such an unlucky teaching situation was uncritically reproduced in (Blomhøj and Kjeldsen 2006): A boy is brought to a hospital with an acute asthmatic attack. There is a drug that can help, if a critical upper concentration in the blood is not surpassed and the concentration doesn't fall below a lower critical level of ineffectiveness. Given the weight of the boy, a dosage scheme might be immediately derived from tables. However, the blood content of two children with the same weight may strongly vary. So the doctor gives one injection and follows the decay of the concentration with 10 measurements in 2-hours intervals over 18 hours. The data table distributed to the pupils shows that the critical lower level is reached already at the third measurement after 6 hours. Certainly with most attentive pupils, I wonder why the teacher's doctor needs 10 measurements to estimate the decay rate of an exponential function when 3 measurements suffice – and let the boy suffer in a perhaps critical state for additional 12 hours.

Personally, I have argued against uninformed naïve modeling patriotism for decades (Booß-Bavnbek 1991, Booß-Bavnbek and Pate 1989a, b). In concrete situations where school administration and teacher qualification only permit the choice between the *pure mathematics-oriented* and the *context-oriented* position one should hesitate to follow too readily the pragmatics of the *context-oriented* or to blame too much the *pure mathematics-oriented* as elitist, outdated and even dead. As in general technology assessment, the conservative approach may be the more humanistic and more soundly futuristic approach also in mathematics education.

## Approach 5 example 3

Unfortunately, the same (otherwise worth reading) article displays a somewhat disorienting ideal scheme of a so-called *mathematical modeling cycle*. By that these authors and thousands of similar articles mean an iterative process of approximating a model to the studied segment of reality in increasing similarity. Contrary to Peirce's realistic quality criterion of modeling as seen from an observer, its goals and its practical success with the model, they lure the reader on a wrong track leading nowhere.

## Approach 5 example 4

Math teachers from all corners of the world have documented their experiences with the contextualization of math teaching. Promising contexts are elements of public health, epidemiology, and capacities of the human brain, exemplified by solutions of mathematical problems or insight in harmonies of numbers, rhythms and accords in music. Other suitable contexts are breeding and crops, harvesting machines, weed and pest control, fertilizers, and the size and value of land property. However, don't let the vitality of chosen contexts take over and remove the abstract concepts to be learned and practiced! And be aware that complex contexts will make ability differences in the student population visible! To address ability differences is not bad but good, but may require additional support in differentiated learning environments.

Often, humanistic and socially engaged teachers select contexts apt for the wanted *empowerment of a critical citizenship* (Freire 1972, Skovsmose 1994), and rightly so. However, there are three fallacies: Firstly, educators please free yourselves from the prejudices of the social sciences regarding mathematical contents. The British Prime Minister Benjamin Disraeli (1804-1881) confirmed that there were two cultures, of the wealthy and the poor in his country, and they still are separated. Later, C.P. Snow (1905-1980) showed that there are similar high fences between sciences and mathematics on one side and the social sciences and humanities on the other. But it is a failure to proliferate that wall between two cultures and deny the natural and intimate relation between class room demands for teaching mathematical content (under the pejorative heading of *numeracy*) vs. class room demands of teaching mathematical meaning (heralded as *mathemacy*). Secondly, educators please minimize your arrogance! Education alone will scarcely transfer power to the underprivileged. Thousands of shy and frightened teachers in the countryside can learn quite a lot from the determination, e.g., of the landless and not vice versa. Thirdly, educators please think constantly about possible optimistic contexts of your teaching instead of silencing your pupils by dwelling too long with negative contexts how realistic they may be. In doubt, teachers will do better when they are perceived of as dreamers rather than as pessimists!

### Approach 5 theoretical explanation
Math teachers have a saying "Don't overestimate the pupils' prior knowledge and don't underestimate their intelligence!" In many places, however, this maxim is double erroneous. Firstly, teachers including relevant contexts will be surprised again and again by the wit and competence of their pupils in their judgments in situations which are familiar to them but remote for the teacher because of their different background. Secondly, it must be repeated once more that our brain is not well suited to handle abstract and formal signs. Humans are not as intelligent as we believe. In the rich kingdom of animals, e.g., only humans – in all historical and ethnological formations - have to be taught the fundamentals of sexual reproduction. Hence, the first criterion for good teaching should be whether the teacher can confirm at the end of an hour that the pupils have a better understanding in some aspects than the teacher had before and that some of the pupils got the possibility to show that they were more clever than their teacher – and not vice versa.

### Approach 5 practical conclusion and claim
Teachers in rural areas have to make the pupils aware of the role of contexts in understanding the content, meaning and value of mathematical concepts, techniques and methods. Contrary to urban life forms with their specific anonymity, perceived mobility and typical exchangeability of relations, rural life forms provide a rich frame for contextualization of the teaching and learning.

## Approach 6 – Of course, at the end only results count. To achieve them the process gains importance. Promises of sample-based tests of learning mathematical experiences.
Addressing the role of tests is a good entrance for discussing goals, methods and experiences in teaching mathematical experiences. There are many reservations against tests. E.g., tests may measure something easily measurable but not very relevant and so become irrelevant. Tests may steer teaching and so become pernicious. That said, tests remind the teacher, the pupils and the parents that there should be a positive result of the teaching act, comparable to putting money in a saving bank. Freire (1972) ridiculed that *banking pedagogy*, and certainly rightly so regarding many classroom situations. However, his imperative to replace result-oriented teaching by *process-oriented teaching* may be valuable for large domains of humanities and social sciences but misses the point of math and science education, where reliable results count. My reservations against the well-intended pupil activity- and inquiry-oriented education are in line with the carefully researched and accommodatingly worded analysis of the necessary balance between product and process orientation in (Sjøberg 2009, 2005, pp. 430-439).

### Approach 6 example 1
The global educational reform thinking includes an assumption that competition, choice, and more-frequent external testing are prerequisites to improve the quality of education. The PISA database shows that those education systems where competition, choice, and test-based accountability have been the main drivers of educational change have not shown progress in international comparisons. Finland has.

The Finns say *test less, learn more*: "Education policies in Finland (in recent years) emphasized teacher professionalism, school-based curriculum, trust-based educational leadership, and school collaboration through networking. Finland has, unlike any other nation, improved its average performance from its already high level in 2000. Although this does not constitute evidence of the failure of test-driven educational reform policies per se, it suggests that frequent standardized student testing is *not* a necessary condition for improving the quality of education." (Sahlberg 2011, p. 66)

### Approach 6 example 2
I share the Finnish idea that testing itself is not a bad thing. I am not an anti-assessment person. Statistically well-designed tests are wonderful means to make teachers attentive to the greater or minor effectiveness of different teaching methods, also relative to different strata of pupils. An in any aspects excellent and professional example is given in (Hansen 2013). Problems arise when the tests become higher in stakes and include sanctions to teachers or schools as a consequence of poor performance. "There are alarming reports from many parts of the world where high-stakes tests have employed as part of accountability policies in education. This evidence suggests that teachers tend to redesign their teaching according to these tests, give higher priority to those subjects that are tested, and adjust teaching methods to drilling and memorizing information rather than understanding knowledge." (Sahlberg 2011, p. 67)

### Approach 6 theoretical explanation
With other subjects like sciences, medicine, and languages, math is a teaching subject where results matter more than process. Process is good, discourse is good. However, when moving forward, pupils will instantly meet abysses where continuation along a desired path is blocked to the uninitiated and untrained. Progress, entering new fields, and climbing to higher levels require a certain security in the delicacies of the mathematical concepts and methods applied. As mentioned before, (Sjøberg 2009, 2005, pp. 430-439) elaborates profoundly on these objections against the well-intended, but unfounded dominance of process orientation in science (and math) teaching.

It is true; there is no *royal route* in mathematics: Gauss gave four radically different proofs for the fundamental theorem of algebra. He excessively used his freedom of inventing fruitful mathematical concepts and to arrange mathematical arguments in a new, surprising, productive and meaningful order. But the goal was always the same: to provide intelligible, convincing, and reliable arguments for the claimed existence of the roots, respectively, possibility of factorization. We as learners and his readers have also plenty of freedom. We can choose the problems, methods and proofs we want. However, in one aspect neither Gauss nor we have any freedom; we have not the choice to believe a little more in the one proof and a bit less in another one; if we doubt one proof, we must check it and either point to a gap or error in the chain of arguments – or confirm the reliability of the arguments. A mathematical argument is valuable, only if it is reliable. Errors can be inspiring, but must be refuted in the very end. There is no place for relativism. Futile arguments definitely have a lesser status than reliable arguments. There is no chance for equality, discourse, and compromise when reliability or unreliability has been tested.

How can we be so sure about mathematical results? Aren't there some logical cracks in the foundations of mathematics, as pointed out in Gödel's Incompleteness Theorem (Manin 2010)? What about the many unsolved mathematical problems? Could they be candidates for true statements that can never be proven? What about the logical status of simple concepts like the set of all subsets of the integers; the real numbers; compact sets; Zorn's lemma? These are interesting questions for the logician. Most mathematicians, however, don't need *foundations* for their work.

They need reliability of their arguments; they need and provide trust in the correctness of results and convince their peers. That is only in seldom cases something a machine can do. The professionalism of mathematicians is that they can say *stop – now the result is reliable*. Also that has to be learned in math classes.

### Approach 6 practical conclusion and claim

The result orientation is mandatory for all math learning. It requires a sharp distinction between the strictly limited periods of the laborious, demanding, sometimes painful acquisition of new mathematical knowledge (depicted in Escher's front cover of Bleecker and Booß-Bavnbek 2013) and the feel-good periods of spare time and leisure. As the Finns say, *teaching less is more.*

Hopefully, in rural environment, the inherent conflict between process- and result-orientation in math education has better chances of a productive solution than in most urban environments. Following the Finnish lessons, it may be easier to implement such a deliberately confined education on the countryside with pupils and parents being more down to earth, than in urban environment with too large an abundance of choices. Pupils will be strengthened by quick response to their advances and failures. All failures should be recognized as intermediate and provisional.

On the basis of such a clear agenda, namely that the main goal of the math classes is learning of contents, of concepts, methods and meaning, not of discourses, the teacher can set the scene for context-oriented math education and pupil participation as the precondition for continued and comprehensive understanding, thought-liberation, intellectual dignity, practical usefulness, and political emancipation.

*Acknowledgements*. I'm indebted to Júlio Faria Corrêa (Roskilde and Campinas) and Philip J. Davis (Brown University, Providence R.I.) for inspiring discussions on the topics of this paper.

### References

Aleksandrov, A.D., Kolmogorov, A.N., and Lavrent'ev, M.A. (1963). *Mathematics: Its content, methods, and meaning*. Vol. I-III. Russian 1956. Translated by K. Hirsch, S. H. Gould, and T. Bartha, omitting though the two sections *The essential nature of mathematics* and *The laws of the development of mathematics* of the original Russian text due to the self-imposed censorship of the translators against "dialectical materialism", Vol I, p. 64 n. Cambridge, Mass.: The M.I.T. Press.

Bleecker, David and Booß-Bavnbek, Bernhelm (2013). *Index theory with applications to mathematics and physics*. Boston: International Press.
http://intlpress.com/site/pub/pages/books/items/00000398/reviews/index.html

Blomhøj, Morten and Kjeldsen, Tinne Hoff (2006). 'Teaching mathematical modelling through project work - Experiences from an in-service course for upper secondary teachers'. *ZDM Mathematics Education* **38**/2, 163-177.

Booß-Bavnbek, Bernhelm (1991): 'Against ill-founded, irresponsible modelling'. In: Niss, M., Blum, W., Huntley, I., (eds.), *Teaching of mathematical modelling and applications* (pp. 70-82). Chichester: Ellis Horwood.

Booß-Bavnbek, Bernhelm and Davis, Philip J. (2013). 'Unity and disunity in mathematics'. *Europ. Math. Soc. Newsletter* **87**/March, 28-31.

Booß-Bavnbek, Bernhelm and Pate, Glen (1989a): 'Expanding risk in technological society through progress in mathematical modeling'. In: C. Keitel, A. Bishop, P. Damerow, P. Gerdes (eds.), *Mathematics, education, and society* (pp. 75-78). Paris: UNESCO, Division of Science, Technical and Environmental Education.

-----, ----- (1989b): 'Information technology and mathematical modelling, the software crisis, risk and educational consequences'. *ZDM Mathematics Education* **21**/5, 167-175.

D'Ambrosio, Ubiratan (2006). *Ethnomathematics. Link between traditions and modernity*. (English translation by Anne Kapple of *Etnomátematica. Elo entre as tradições e modernidade,* 2002, Coleção Tendências em Educação Matemática, Editoria Autêntica: Belo Horizonte). Rotterdam: Sense Publishers. xi+89 pages.

Damerow, Peter and Wolfgang Lefèvre (1981). *Rechenstein, Experiment, Sprache. Historische Fallstudien zur Entstehung der exakten Wissenschaften*. Stuttgart: Klett-Cotta.

Davis, Philip J. and Hersh, Reuben (1981). *The mathematical experience*. With an introduction by Gian-Carlo Rota. Boston, Mass.: Birkhäuser.

Freire, Paulo (1972). *Pedagogy of the oppressed*. Translated by Myra Bergman Ramos (Pedagogia do oprimido, manuscrito em português de 1968). New York: Herder and Herder.

Freudenthal, Hans (2002). *Revisiting mathematics education: China lectures.* Dordrecht: Kluwer Academic. First edition 1991. Also: http://p4mriunismuh.files.wordpress.com/2010/08/revisiting-mathematics-education.pdf, accessed Dec. 6, 2013.

Gerdes, Paulus (2000*). Le cercle et le carré: créativité géométrique, artistique et symbolique de vannières et vanniers d'Afrique, d'Amérique, d'Asie et d'Océanie*. Préf. de Maurice Bazin. Paris: L'Harmattan.

Gramsci, Antonio (1971). *Selections from the Prison Notebooks*. New York: International Publishers.

Hansen, Ali (2013). 'Comprehension of text. Can additional verbal explanations enhance students' performance in solving text-based mathematical assignments?' (in Danish). *MONA - Matematik- og Naturfagsdidaktik - Tidsskrift for undervisere, formidlere og forskere* (submitted).

Haug, Wolfgang Fritz (1996). *Philosophieren mit Brecht und Gramsci*. Berlin & Hamburg: Argument-Verlag.

Høyrup, Jens (2002). *Lengths, widths, surfaces: a portrait of Old Babylonian algebra and its kin*. (Studies and Sources in the History of Mathematics and Physical Sciences). New York: Springer.

Husserl, Edmund (1970). *The origin of geometry*. Appendix VI of *The crisis of European sciences and transcendental phenomenology. An introduction to phenomenological philosophy* (pp. 353-378). Evanston, IL: Northwestern University Press. English translation by David Carr of the German original of 1936.

Keitel, Christine, Kotzmann, Ernst, & Skovsmose, Ole (1993). 'Beyond the tunnel vision: Analysing the relationship between mathematics, society, and technology'. In: C. Keitel & K. Ruthven (Eds.), *Learning from computers: Mathematics education and technology* (pp. 243-279). Berlin: Springer.

Kierkegaard, Søren (1959). *Either/Or*. (Translated by David F. Swenson and Lillian Marvin Swenson from the Danish Enten - Eller. Et Livs-Fragment, 1, 1843). Volume I. Princeton: Princeton University Press.

Knijnik, Gelsa and Wanderer, Fernanda (2010). 'Mathematics education and differential inclusion: a study about two Brazilian time–space forms of life'. *ZDM Mathematics Education* **42**, 349–360.

Krantz, Steven G. (1997). *A primer of mathematical writing*. Providence, Rhode Island: American Mathematical Society.

Lenin, Vladimir I. (1909, 1999). *Materialism and empirio-criticism*. Marxists Internet Archive, www.marxists.org/ (Source: *Lenin Collected Works*, Progress Publishers, 1972, Moscow, Vol. **14**, pp. 17-362).

Luria, Alexander Romanovich (1972). *The man with a shattered world*. New York: Basic Books, Inc. (translated by L. Solotaroff).

Manin, Yuri (2010). *A course in mathematical logic for mathematicians*. Second edition. Chapters I–VIII translated from the Russian by Neal Koblitz. With new chapters by Boris Zilber and the author. Graduate Texts in Mathematics, 53. New York: Springer.

Marx, Karl (1844). Introduction to the critique of Hegel's philosophy of right. http://www.marxists.org/archive/marx/works/1843/critique-hpr/intro.htm

-------- (1845). *Theses on Feuerbach*. http://www.marxists.org/archive/marx/works/1845/theses/

Mesquita, Monica (2013). 'Communitarian mathematics education. Walking into boundaries'. *EM TEIA - Revista de Educação Matemática e Tecnológica Iberoamericana* **4**/4, Thematic edition in "Rural education: contributions from mathematics and technological education", this issue.

Miguel, Antonio and Mendes, Iran Abreu (2010). 'Mobilizing histories in mathematics teacher education: memories, social practices, and discursive games'. *ZDM Mathematics Education* **42**/3-4, 381–392.

Minsky, Marvin (1967). *Computation: finite and infinite machines*. Englewood Cliffs, N.J.: Prentice-Hall.

Moon, Bob (1986). *The new maths curriculum controversy: an international story.* Studies in curriculum history series, vol. 5. London: Falmer Press.

Naur, Peter (1992). *Computing: a human activity*. New York: ACM Press/Addison-Wesley.

Otte, Michael (1994). *Das Formale, das Soziale und das Subjektive: Eine Einführung in die Philosophie und Didaktik der Mathematik*. Frankfurt am Main: Suhrkamp Taschenbuch Wissenschaft,

----- (1997). 'Mathematik und Verallgemeinerung - Peirce' semiotisch-pragmatische Sicht'. *Philosophia Naturalis* **34** (2), 175–222.

Peirce, Charles Sanders (1992, 1998). *The Essential Peirce, Selected philosophical writings*. Volume 1 (1867–1893), Nathan Houser and Christian J. W. Kloesel, eds.; Volume 2 (1893–1913), Peirce Edition Project. Bloomington and Indianapolis, IN: Indiana University Press.

------- (1877). 'The Fixation of Belief'. Popular Science Monthly 12 (November 1877), 1-15. http://www.peirce.org/writings/p107.html, light edited on http://milne.ruc.dk/~booss/Mamocalc/MaMoTeach_bbb/Peirce_on_the_Philosophy_of_Modeling.pdf

Pinker, Steven (1994). *The language instinct*. New York: William Morrow.

Radford, Luis (2006). 'The anthropology of meaning'. *Educational Studies in Mathematics* **61**, 39–65.

Rosser, J. Barkley (1953). *Logic for mathematicians*. New York, Toronto, London: McGraw-Hill Book Company Inc. Digitized by the Internet Archive in 2010 at http://archive.org/stream/logicformathemat00ross/logicformathemat00ross_djvu.txt, accessed Dec. 6, 2013.

Sahlberg, Pasi (2011). *Finnish lessons: what can the world learn from educational change in Finland?* Foreword by Andy Hargreaves. The series on school reform. New York: Teachers College Press.

Sánchez, Mario (2010). *How to stimulate rich interactions and reflections in online mathematics teacher education*? PhD-thesis. Roskilde: IMFUFA-tekst 472.

Schubring, Gert (2007). 'Der Aufbruch zum "funktionalen Denken": Geschichte des Mathematikunterrichts im Kaiserreich', *N.T.M.* **15**, 1-17.

Science Obituaries: 'Israel Gelfand'. *The Telegraph*. 26 Oct 2009.

Sfard, Anna (1991). 'On the dual nature of mathematical conceptions: reflections on processes and objects as different sides of the same coin'. *Educational Studies in Mathematics* **22**, 1-36.

Sjøberg, Svein (2009, 2005). *Naturfag som allmenndannelse: en kritisk fagdidaktikk* [Norwegian, *Science education as general knowledge (German: Allgemeinbildung): a critical subject-related* *education (German: Fachdidaktik)*]. 3d edition. Oslo: Gyldendal akademisk. Quoted from the Danish translation by Alf Andersen, 2005, Aarhus: Klim.

Skovsmose, Ole (1994). *Towards a philosophy of critical mathematics education*. Dordrecht: Kluwer Academic Publishers.

Wittgenstein, Ludwig (1953/2009). *Philosophical investigations*. Includes the original German text in addition to the English translation by G.E.M. Anscombe, P.M.S. Hacker and J. Schulte. Oxford: Blackwell Publishing.